\documentstyle{article}
\begin{document}
\LARGE
\begin{center}
A Note on a Specific Pencil of Conics in the Galois Fields of Order
     2$^{n}$
\end{center}
\vspace*{-0.1cm}
\Large
\begin{center}
Metod Saniga\\

\vspace*{.2cm}
\normalsize
Astronomical Institute,  Slovak Academy of Sciences,
SK-059 60 Tatransk\' a Lomnica,\\ The Slovak Republic\\
{\it E-mail: msaniga@astro.sk}
\end{center}

\vspace*{.1cm}
\normalsize
\noindent
{\bf Abstract}\\ 

The pencil of conics featuring
three degenerate conics each of which is a line-pair is briefly inspected in 
a Galois field of characteristic two. It is shown that if two degenerates
are conjugate imaginary line pairs, the third must be a {\it real} line pair;
this contradicts Campbell's claim (Campbell 1927) that all the three singular
conics are conjugate imaginary line pairs.

\vspace*{0.8cm}
\large
\noindent
Our remark concerns the structure of the pencil of conics containing three
degenerate conics of which two represent conjugate imaginary line pairs,
i.e. the pencil defined by Eq. (24) in [1]. In what follows we
will demonstrate that Campbell's claim that also the third singular conic  
`must be a conjugate imaginary line pair' is false, because this conic is,
in fact, a {\it real} line pair.

The proof relies on the following theorem: the expression
\begin{equation}
u^{2} + v^{2} + \Theta uv,~~~~~~\Theta \neq 0,
\end{equation}
where $u$ and $v$ are regarded as variables and $\Theta$ is a parameter, is
reducible or irreducible in the Galois field of order $2^{n}$ (GF(2$^{n}$)
iff, respectively, 
\begin{equation}
D \left(1/ \Theta^{2} \right) = 0,
\end{equation}
or
\begin{equation}
D \left(1/ \Theta^{2} \right) = 1,
\end{equation}
where
\begin{equation}
D(w) \equiv w + w^{2} + w^{4} + \ldots + w^{2^{n-1}}
\end{equation}
 (see, e.g. [2]). The pencil of conics concerned
is (see Eq. (24) of [1])
\begin{equation}
C(\lambda, \mu) \equiv \lambda C_{1} + \mu C_{2} \equiv
\lambda \left( x^{2} + y^{2} + \alpha xy \right) +
\mu \left( x^{2} + z^{2} + \beta xz \right),
\end{equation}
where
\begin{equation}
\alpha \beta \neq 0,~~~\alpha \neq \beta.
\end{equation}
If the conics $C_{1} = 0$ and $C_{2} = 0$ are to represent conjugate imaginary
line pairs, both $C_{1}$ and $C_{2}$ must be {\it ir}reducible, i.e.
\begin{equation}
D\left(1/ \alpha^{2} \right) = 1 = D\left(1/ \beta^{2} \right).
\end{equation}
Our task is to find the character of the third degenerate conic of the pencil,
which is given by (see p. 405 of [1])
\begin{equation}
C_{3} \equiv \left( \beta^{2} + \alpha^{2} \right) x^{2} + \beta^{2} y^{2} +
\alpha^{2} z^{2} + \beta^{2} \alpha x y + \beta \alpha^{2} x z = 0.
\end{equation}
To this end in view we first notice that with the help of the relation
\begin{equation}
(u + v)^{2} = u^{2} + v^{2}
\end{equation}
Eq. (8) can be cast into the form
\begin{equation}
\bar{C}_{3} = x^{2} + \left( \frac{\beta y + \alpha z}{\alpha + \beta}
\right)^{2} + \frac{\alpha \beta}{\alpha + \beta}~ x~  \frac{\beta y + 
\alpha z}{\alpha + \beta} = 0,
\end{equation}
which, using a non-singular transformation
\begin{eqnarray}
x' &=& x, \nonumber \\ \nonumber \\
y' &=& \frac{\beta y + \alpha z}{\alpha + \beta}, \nonumber \\ \nonumber \\
z' &=& \frac{\alpha y + \beta z}{\alpha + \beta},
\end{eqnarray}
is sent into
\begin{equation}
C'_{3} = x'^{2} + y'^{2} + \gamma x' y' = 0,
\end{equation}
where
\begin{equation}
\gamma \equiv \frac{\alpha \beta}{\alpha + \beta}.
\end{equation}
The shape of $C'_{3}$ is identical with that of Eq. (1) so that its
character, in the light of the above-introduced theorem, depends solely on 
the value of 
$D\left( 1/ \gamma^{2} \right)$. In order to find the latter we first observe 
that
\begin{equation}
D\left( 1/ \gamma^{2} \right) = D  \left( \frac{ \alpha^{2} + \beta^{2}}
{\alpha^{2} \beta^{2}} \right) = D \left( \frac{1}{\alpha^{2}} +
\frac{1}{\beta^{2}} \right).
\end{equation}
Further, from the definition of $D(w)$ and Eq. (9) it can easily be verified 
that
\begin{equation}
D(u + v) = D(u) + D(v),
\end{equation}
which implies that
\begin{equation}
D\left( 1/ \gamma^{2} \right) = D \left( 1/\alpha^{2} \right) +
D \left( 1/\beta^{2} \right) = 1 + 1 = 0
\end{equation}
where we also took into account Eq. (7) and the fact that $w + w$ = 0 for 
any $w \in $ GF(2$^{n}$). Eq. (16) tells us that $C'_{3}$ is
{\it reducible} and the corresponding conic $C'_{3} = 0$ thus, indeed, 
represents a {\it real} line pair.

\end{document}